\documentclass[a4paper,12pt,leqno]{article}
\usepackage{latexsym}

\usepackage{amssymb} 
\usepackage{amsmath} 
\usepackage{theorem}

\def\Z{{\mathbb{Z}}}
\def\K{{\mathbb{K}}}
\def\C{{\mathbb{C}}}
\def\A{{\mathcal{A}}}

\DeclareMathOperator{\rank}{rank}
\DeclareMathOperator{\codim}{codim}

\DeclareMathOperator{\Der}{Der}

\DeclareMathOperator{\res}{res}

\def\d{{\mathbf{d}}}
\def\dhA{{\d_{H_0}\A}}
\def\Hzero{{H_0}}
\DeclareMathOperator{\pdim}{pdim}

\numberwithin{equation}{section}

\newcommand{\owari}{\hfill$\square$}

\theoremstyle{break}
\newtheorem{theorem}{Theorem}[section]
\newtheorem{prop}[theorem]{Proposition}
\newtheorem{cor}[theorem]{Corollary}

\newtheorem{define}[theorem]{Definition}
\newtheorem{rem}[theorem]{Remark}

\newtheorem{example}[theorem]{Example}

\title{Free arrangements and coefficients of characteristic polynomials}
\author{Takuro Abe\thanks{Department of Mechanical Engineering and Science, 
Kyoto University, Yoshida Honmachi, Sakyo-ku, Kyoto 6068501, Japan. 
email: abe.takuro.4c@kyoto-u.ac.jp} and 
Masahiko Yoshinaga\thanks{Department of Mathematics, 
Kyoto University, 
Kitashirakawa Oiwakecho, Sakyo-ku, Kyoto, 606-8502, Japan. 
email: mhyo@math.kyoto-u.ac.jp
}}
\date{\today} 

\begin{document}

\maketitle

\begin{abstract}
Ziegler showed that free arrangements 
have free restricted multiarrangements (multirestrictions). 
After Ziegler's work, 
several results concerning ``reverse direction'', namely 
characterizing freeness of an arrangement via that of 
multirestriction, have 
appeared. In this paper, we prove that the second Betti 
number of the arrangement plays a crucial role. 
\end{abstract}

\section{Introduction}
\label{sec:intro}

Let $V$ be a vector space of dimension $\ell$ over 
a field $\K$. Fix a system of coordinate 
$(z_1, \dots, z_\ell)$ of $V^*$. We denote by 
$S=S(V^*)=\K[z_1, \dots, z_\ell]$ the 
symmetric algebra. 
A hyperplane arrangement $\A=\{H_1, \dots, H_n\}$ is 
a finite collection of hyperplanes in $V$. 

Freeness of an arrangement is a key notion which 
connects arrangement theory with 
algebraic geometry and combinatorics. There are several 
ways to prove freeness, e. g. using Saito's criterion 
\cite{Sa2}, addition-deletion theorem \cite{ter-arr}, etc. 
In \cite{zie-multi}, Ziegler proved that the multirestriction 
$(\A^{H_0}, m^{H_0})$ of a free arrangement $\A$ is also free 
(see \S\ref{sec:chi} for details). 
The converse is not true in general. 
However Schulze \cite{Sc} recently proved that if the dimension 
is $\ell\leq 4$ (or $\ell\geq 5$ under tameness assumption), 
freeness of $\A$ is characterized in terms of 
multirestriction and characteristic 
polynomials. The purpose of this paper is to give a 
stronger characterization of freeness for any dimension. 
Namely, we characterize the freeness in terms of the multirestriction 
and the second coefficients of characteristic polynomials 
(without posing any conditions on 
dimension or tameness). 

This paper is organized as follows. In \S\ref{sec:chi}, 
we recall basic facts on characteristic polynomials for 
both simple and multiarrangements. In \S\ref{sec:viamulti}, 
we recall results from \cite{Y1, Y2, Sc}, which will be 
used in the proof of the main result. 
In \S\ref{sec:restr} we formulate a new combinatorial 
technique. Fix a hyperplane $H_0\in\A$. Then we can associate 
two arrangements: deconing $\dhA$ and the restriction $\A^{H_0}$. 
We define a natural map $\rho: L(\dhA)\rightarrow 
L(\A^{H_0})$ of their intersection posets. The map 
will be used in the proof of main result. 
In \S\ref{sec:main}, we state and prove the main result. 
In \S\ref{sec:related}, we prove several related results 
by localizing our main result.

\medskip

\noindent
{\bf Acknowledgement} Parts of this work were done 
during the Workshop on Free Divisors 
(Warwick, June 2011) and 
Hyperplane arrangements and applications 
(Vancouver, August 2011). 
The authors thank organizers and participants of 
these conferences. Especially the similarity 
with formality (Remark \ref{rem:formal}) was pointed 
out by Michael Falk. 
Both authors are supported by 
JSPS Grant-in-Aid for Young Scientists (B). 

\section{Characteristic polynomials}
\label{sec:chi}

In this section, we summarize several facts on 
the characteristic polynomials. 

\subsection{For simple arrangements}
\label{sec:simple}

Let $\A$ be an arrangement of affine hyperplanes 
in an affine space $V$ of $\dim_\K V=\ell$. 
Let $L(\A)$ be the set of nonempty 
intersections of elements of $\A$. Define a 
partial order on $L(\A)$ by $X\leq Y\Longleftrightarrow 
X\supseteq Y$, where $X, Y\in L(\A)$. 
Then $L(\A)$ is a ranked poset with 
$\rank X=\codim X$. 
Denote the set of $X\in L(\A)$ of rank $r$ by 
$L_r(\A)=\{X\in L(\A)\mid \codim X=r\}$. 

Let $\mu:L(\A)\rightarrow\Z$ be the M\"obius function 
of $L(\A)$ defined by $\mu(V)=1$, and for 
$X>V$ by the recursion $\sum_{Y\leq X}\mu(Y)=0$. 
The characteristic polynomial of $\A$ is defined as 
$\chi(\A, t)=\sum_{X\in L(\A)}\mu(X)t^{\dim X}\in\Z[t]$. 
Set the coefficients of the characteristic polynomial 
$b_k(\A)$ as 
$\chi(\A, t)=t^{\ell}-b_1(\A)t^{\ell-1}+b_2(\A)t^{\ell-2}-
\dots+(-1)^{\ell}b_{\ell}(\A)$. 
The following Local-Global formula is straightforward. 

\begin{prop}[Local-Global formula for Betti numbers]
\label{prop:betti}
Let $\A$ be an affine arrangement. Then 
\begin{equation}
\label{eq:betti}
b_k(\A)=
\sum_{X\in L_k(\A)}b_k(\A_X), 
\end{equation}
where $\A_X=\{H\in\A\mid H\supset X\}$. 
\end{prop}

We call $\A$ central if $\bigcap_{H\in\A}H\neq\emptyset$. 
In this case, the characteristic 
polynomial $\chi(\A, t)$ is divisible by $(t-1)$. 
We denote $\chi_0(\A, t)=\frac{1}{t-1}\chi(\A, t)$. 

Given a nonempty central arrangement $\A$ and 
a hyperplane $H_0\in\A$. Choose coordinates 
$z_1, \dots, z_\ell$ of $V$ so that 
$H_0=\{z_\ell=0\}$. Let $H_0'=\{z_\ell=1\}$ be 
an affine hyperplane parallel to $H_0$. 
The deconing $\d_{H_0}\A$ of $\A$ with respect to $H_0$ is 
the affine arrangement $H_0'\cap\A$ on $H_0'$. 
The deconing $\d_{H_0}\A$ is an affine arrangement of 
rank $(\ell-1)$ whose characteristic polynomial 
satisfies 
$\chi(\d_H\A, t)=\chi_0(\A, t)$.

\subsection{For multiarrangements}
\label{subsec:multi}


Let $\A$ be a central arrangement. 
A map $m:\A\rightarrow\Z_{\geq 0}$ is called a 
multiplicity. 
We define the $S$-module $D(\A, m)$ for a multiarrangement 
$(\A, m)$ by 
$$
D(\A, m)=
\{
\delta\in\Der_\K(S)\mid
\delta(\alpha_H)\in(\alpha_H^{m(H)}), 
\forall H\in\A
\}, 
$$
where 
$\alpha_H$ is a linear form such that 
$\ker \alpha_H=H$ for each hyperplane $H\in\A$. 
A multiarrangement $(\A, m)$ is called free 
with exponents $(d_1, \dots, d_\ell)$ 
if $D(\A, m)$ is a free $S$-module with a homogeneous 
basis $\delta_1, \dots, \delta_\ell\in D(\A, m)$ 
such that $\deg\delta_i=d_i$.

Let $\Omega^1_V=\bigoplus_{i=1}^\ell S\cdot dx_i$ and 
$\Omega_V^p=\bigwedge^p\Omega^1_V$. 
We define an $S$-module $\Omega^p(\A, m)$ of a multiarrangement 
$(\A, m)$ by 
$$
\Omega^p(\A, m)=\left\{
\omega\in
\frac{1}{Q}\Omega_V^p\mid 
\frac{d\alpha_H}{\alpha_H^{m(H)}}\wedge\omega\in
\frac{1}{Q}\Omega_V^{p+1}
\right\}, 
$$
where 
$Q=Q(\A, m)=\prod_{H\in\A}\alpha_H^{m(H)}$. 
Next we recall the characteristic polynomial of 
a multiarrangement $(\A, m)$ \cite{ATW}. 
Recall that for a finitely generated graded 
$S$-module $M=\bigoplus_{d\in\Z}M_d$, 
the Hilbert series $P(M, x)\in\Z[x^{-1}][[x]]$
of $M$ is defined by 
$$
P(M, x)=
\sum_{d\in\Z}(\dim_\K M_d)x^d. 
$$
For a multiarrangement $(\A, m)$, we define 
$$
\Phi((\A, m); x, t)=
\sum_{p=0}^\ell P(\Omega^p(\A, m), x)(t(1-x)-1)^p
\in\Z[x^{\pm 1}, t]. 
$$

\begin{define}[\cite{ATW}]
\label{def:atw}
The characteristic polynomial of $(\A, m)$ 
is defined as follows. 
$$
\chi((\A, m), t)=
\lim_{x\rightarrow 1}\Phi((\A, m); x, t)\in\Z[t]. 
$$
\end{define}
Define the integer $\sigma_i(\A, m)\in\Z$ by 
$$
\chi((\A, m), t)=
t^\ell
-\sigma_1(\A, m)t^{\ell-1}
+\sigma_2(\A, m)t^{\ell-2}
-
\dots
+(-1)^{\ell}\sigma_{\ell}(\A, m). 
$$

\begin{rem}
To be precise, the characteristic polynomial 
$\chi((\A, m), t)$ was first defined by using 
the dual module $D^p(\A, m)$ in \cite{ATW}. Set 
$$
\psi((\A, m), t, q)=
\sum_{p=0}^\ell
P(D^p(\A, m), q)(t(q-1)-1)^p, 
$$
and it was defined that 
$\chi((\A, m), t):=(-1)^\ell
\lim_{q\rightarrow 1}\psi((\A, m), t, q)$. 
However, this is equivalent to Definition 
\ref{def:atw}. 
It is proved by checking the following two facts. 
\begin{itemize}
\item[(i)] 
``Local-global formula'' (Prop. \ref{prop:locgl} below) holds 
for both $\Phi((\A, m); 1, t)$ and $\psi((\A, m); 1, t)$. 
\item[(ii)] 
The constant terms coincide. 
\end{itemize}
The assertion (i) can be proved in a similar way with 
\cite{ATW}. 
To verify the assertion (ii), we have to prove 
that 
$\Phi((\A, m); 1, 0)=\sum_{p=0}^\ell (-1)^p 
P(\Omega^p(\A, m), x)$ and 
$\psi((\A, m); 1, 0)=\sum_{p=0}^\ell (-1)^{\ell-p} 
P(D^p(\A, m), x)$ are equal, which is proved by using 
the isomorphism 
$
D^p(\A, m)\stackrel{\simeq}{\longrightarrow}
\Omega^{\ell-p}(\A, m)[-\deg Q]: \delta\longmapsto
\iota_\delta\left(
\frac{dx_1\wedge\dots\wedge dx_\ell}{Q}\right)
$. 
By (i), every coefficient of $\chi((\A, m), t)$ can 
be interpreted as the sum constant terms of certain 
localized subarrangements. The constant terms of these 
localizations coincides thanks to (ii). 
\end{rem}

\begin{rem}
It is not known whether or not 
$\sigma_i(\A, m)\geq 0$ holds. 
In \cite{A2}, it is proved under the assumption of 
tameness. 
\end{rem}

Let $X\in L(\A)$. Recall that 
$\A_X$ is the set of all hyperplanes which contains $X$. 
The restricted multiplicity is denoted by 
$m_X=m|_{\A_X}$. 

\begin{prop}[\cite{ATW}]
\label{prop:locgl}
Let $(\A, m)$ be a central multiarrangement. Then 
\begin{equation}
\label{eq:sigma}
\sigma_k(\A, m)=
\sum_{X\in L_k(\A)}
\sigma_k(\A_X, m_X). 
\end{equation}
\end{prop}

Let $X\in L_2(\A)$. Since $(\A_X, m_X)$ is 
a rank two multiarrangement, it is free. 
We denote the exponents by 
$(e_1(X), e_2(X), 0, \dots, 0)$. 
Using these data, $\sigma_1(\A, m)$ and $\sigma_2(\A, m)$ 
are expressed as follows. 
\begin{equation}
\label{eq:lmp2}
\begin{split}
\sigma_1(\A, m)&=\sum_{H\in\A}m(H), \\
\sigma_2(\A, m)&=\sum_{X\in L_2(\A)}e_1(X)\cdot e_2(X). 
\end{split}
\end{equation}

If $(\A, m)$ is free with exponents 
$(d_1, \dots, d_\ell)$, then 
$$
\chi((\A, m), t)=\prod_{i=1}^\ell(t-d_i). 
$$

\section{Freeness via multirestriction}
\label{sec:viamulti}

Multiarrangements naturally appear as 
restriction of simple arrangements. Let $\A$ be a 
central arrangement. Fix a hyperplane $H_0\in\A$. Let 
$Q(\A')$ be the defining equation 
of the deleted arrangement $\A'=\A\setminus\{H_0\}$. 
Ziegler's multirestriction is a multiarrangement on $H_0$ 
defined by the equation $Q(\A')|_{H_0}$. We denote 
it by $(\A^{H_0}, m^{H_0})$. We have 
$Q(\A^{H_0}, m^{H_0})=Q(\A')|_{H_0}$. 
\begin{prop}[\cite{zie-multi}]
\label{prop:ziegler}
Let $\A$ be a free arrangement. 
Let $H_0\in\A$ and let $(\A^{H_0}, m^{H_0})$ be the 
Ziegler restriction. Then 
\begin{equation}
\label{eq:equal}
\chi_0(\A, t)=
\chi((\A^{H_0}, m^{H_0}), t). 
\end{equation}
In other words, 
\begin{equation}
\label{eq:beqs}
b_k(\dhA)=\sigma_k(\A^{H_0}, m^{H_0}), 
\end{equation}
for $k=1, \dots, \ell-1$. 
\end{prop}
In general, (\ref{eq:equal}) does not hold. 
However to characterize the freeness, 
the formula (\ref{eq:equal}) plays an important role as follows. 

\begin{prop}[\cite{Y2}]
\label{prop:3arr}
Let $\ell=3$. 
Let $\A$ be a central arrangement in $\K^3$ and 
$H_0\in\A$. 
\begin{itemize}
\item[(1)] $b_2(\dhA)\geq \sigma_2(\A^{H_0}, m^{H_0})$. 
\item[(2)] $\A$ is free if and only if 
$b_2(\dhA)=\sigma_2(\A^{H_0}, m^{H_0})$. 
\end{itemize}
\end{prop}

\noindent
When $\ell=3$, the 
condition $b_2(\dhA)=\sigma_2(\A^{H_0}, m^{H_0})$ 
is equivalent to the coincidence of 
characteristic polynomials: 
\begin{equation}
\label{eq:eqchi}
\chi(\dhA, t)=\chi((\A^{H_0},m^{H_0}),t). 
\end{equation}
The result above has been generalized by 
M. Schulze as follows. 

\begin{prop}[\cite{Sc}]
\label{prop:sc}
Let $\A$ be a central arrangement in $\K^\ell$ and 
fix $H_0\in\A$. 
Suppose that $\ell=4$ or 
$\ell\geq 5$ with (weakly) tameness assumption (i.e., 
$\pdim\Omega^p(\A^{H_0}, m^{H_0})\leq p,\ \forall p$). 
Then $\A$ is free if and only if 
\begin{itemize}
\item $(\A^{H_0}, m^{H_0})$ is free, and 
\item the relation (\ref{eq:eqchi}) holds. 
\end{itemize}
\end{prop}

\begin{define}
\label{def:along}
$H_0\in\A$. $\A$ is said to be {\em locally free along 
$H_0$} if $\A_X$ is free for all $X\in L(\A)$ with 
$0\neq X\subset H_0$. 
\end{define}

\begin{prop}[\cite{Y1}]
\label{prop:4arr}
Let $\A$ be a central arrangement in $\K^\ell$ and 
fix $H_0\in\A$ with $\ell\geq 4$. 
Then $\A$ is free if and only if 
\begin{itemize}
\item $(\A^{H_0}, m^{H_0})$ is free, and 
\item $\A$ is locally free along $H_0$. 
\end{itemize}
\end{prop}
In \S\ref{sec:main}, we generalize Proposition 
\ref{prop:3arr} and \ref{prop:sc} to higher dimensions. 

\section{Combinatorial restriction map}
\label{sec:restr}

Let $\A$ be a central arrangement in $V=\K^\ell$. 
Fix $H_0\in\A$. 
Recall that $\d_{H_0}\A$ is an affine arrangement in 
$H_0'$. Hence we may consider $X\in L_k(\d_{H_0}\A)$ to be 
an affine subspace of $V$ of dimension $(\ell-k-1)$. 
(Note that by definition, $X\in L_k(\d_{H_0}\A)$ is codimension 
$k$ in $H_0'$.) 
Then $X$ generates a $(\ell-k)$-dimensional linear subspace 
$\K X\subset V$. By taking intersection with 
$H_0$, we obtain an $(\ell-k-1)$-dimensional 
linear subspace $\K X\cap H_0\in L_k(\A^{H_0})$ of $H_0$. 
We denote the map by $\rho$
\begin{equation*}
\rho:L(\d_{H_0}\A)\longrightarrow
L(\A^{H_0}):\ 
X\longmapsto \K X\cap H_0
\end{equation*}
which preserves the rank and the order of posets. 

The map $\rho$ is compatible with 
the localization in the following manner. 
Let $X\in L(\A)$ with $X\subset H_0$. Then 
$H_0\in\A_X\subset\A$, and the following diagram commutes. 
\begin{equation}
\begin{array}{ccc}
L(\dhA)&\stackrel{\rho}{\longrightarrow}&L(\A^{H_0})\\
\uparrow&&\uparrow\\
L(\dhA_X)&\stackrel{\rho_X}{\longrightarrow}&L(\A_X^{H_0}).
\end{array}
\end{equation}

Furthermore the vertical maps are full, that is, 
$Y_1\in L(\dhA_X), Y_2\in L(\dhA)$ with 
$Y_2\leq Y_1$, then $Y_2\in L(\dhA_X)$. 
This implies that if $Y_1\in L(\dhA_X)$, then 
the value M\"obius function of $Y_1$ in 
$L(\dhA_X)$ is equal to that in 
$L(\dhA)$.

\section{Main results}
\label{sec:main}

\begin{theorem}
\label{thm:main}
Let $\A$ be a central arrangement in $V=\K^\ell$ 
(with $\ell\geq 3$) and 
$H_0\in\A$. 
\begin{itemize}
\item[(1)] $b_2(\d_{H_0}\A)\geq\sigma_2(\A^{H_0}, m^{H_0})$. 
\item[(2)] The equality 
$b_2(\d_{H_0}\A)=\sigma_2(\A^{H_0}, m^{H_0})$ holds if 
and only if 
$\A_X$ is free for all $X\in L_3(\A)$ with $X\subset H_0$. 
(We may say that $\A$ is locally free in codimension three along 
$H_0$.)
\item[(3)] 
Assume that $(\A^{H_0}, m^{H_0})$ is free. Then 
the following are equivalent. 
\begin{itemize}
\item[(i)] $\A$ is free. 
\item[(ii)] $\chi(\dhA, t)=\chi((\A^{H_0}, m^{H_0}), t)$. 
\item[(iii)] $b_2(\d_{H_0}\A)=\sigma_2(\A^{H_0}, m^{H_0})$. 
\item[(iv)] $\A$ is locally free in codimension three along $H_0$. 
\end{itemize}
\end{itemize}
\end{theorem}

{\em Proof.} We first prove (1) and (2). 
By using 
$L_2(\dhA)=\bigsqcup\limits_{X\in L_2(\A^{H_0})}\rho^{-1}(X)$, 
we have 
\begin{equation}
\label{eq:dec}
\begin{split}
b_2(\dhA)&=\sum_{Y\in L_2(\dhA)}\mu(Y)\\
&=
\sum_{X\in L_2(\A^{H_0})}\left(\sum_{Y\in\rho^{-1}(X)}\mu(Y)\right)\\
&=\sum_{X\in L_2(\A^{H_0})}b_2(\d_{H_0}(\A_X)). 
\end{split}
\end{equation}
Suppose $X\in L_2(\A^{H_0})$. 
Then 
$X\in L_3(\A)$ with $X\subset H_0$, and 
$\A_X$ is a rank $3$ central arrangement. 
Hence from Proposition \ref{prop:3arr}, 
we have 
\begin{equation}
\label{eq:locineq}
b_2(\d_{H_0}(\A_X))\geq e_1(X)\cdot e_2(X), 
\end{equation}
where $(e_1(X), e_2(X))$ is the exponents of 
$(\A^{H_0}_X, m^{H_0}_X)$ as in (\ref{eq:lmp2}). 
Hence we have 
\begin{equation}
\label{eq:locineq2}
b_2(\dhA)\geq\sum_{X\in L_2(\A^{H_0})} e_1(X)\cdot e_2(X)
=\sigma_2(\A^{H_0}, m^{H_0}). 
\end{equation}
Thus (1) is proved. Furthermore, in 
(\ref{eq:locineq}), the equality holds if and only if 
$\A_X$ is free. 
Hence the equality 
$b_2(\d_{H_0}\A)=\sigma_2(\A^{H_0}, m^{H_0})$ holds 
if and only if $\A_X$ is free for all $X\in L_2(\A^{H_0})$. 
Thus we have (2). 

Next we prove (3). 
First, (i)$\Longrightarrow$(ii)$\Longrightarrow$(iii)
$\Longleftrightarrow$(iv) is obvious (from (2)). 
We shall prove (iv)$\Longrightarrow$(i) by induction 
on $\ell$. 
By definition, 
$\A_X$ is free for all $X\in L_3(\A)$ with $X\subset H_0$. 
If $\ell=3$, then $L_3(\A)=\{0\}$, $\A_X=\A$ for $X\in L_3(\A)$ 
and there is nothing to prove. 
Let $\ell\geq 4$. We will use 
Proposition \ref{prop:4arr}. 
It suffices to show that $\A$ is locally free along $H_0$. 
Let $Z\in L(\A)$ with $0\neq Z\subset H_0$. 
Then $\A_Z$ has rank at most $(\ell-1)$ since $Z\neq 0$. 
Since $(\A_Z^{H_0}, m_Z^{H_0})$ is a localization of 
$(\A^{H_0}, m^{H_0})$, it is free. It is easily checked that 
$\A_Z$ satisfies (iii). 
Hence by the inductive assumption, $\A_Z$ is free. 
Consequently, $\A$ is locally free along $H_0$. 
\owari

The following are immediate. 

\begin{cor}
\label{cor:main}
Let $\A$ and $H_0\in\A$ be as above. Suppose that 
$(\A^{H_0}, m^{H_0})$ is free with exponents $(d_1, \dots, d_{\ell-1})$. 
Then the inequality 
$$
b_2(\dhA)\geq
\sum_{1\leq i<j\leq\ell-1}d_id_j 
$$
holds. Furthermore, $\A$ is free if and only if the equality holds. 
\end{cor}

\begin{cor}
Let $\A_1$ and $\A_2$ be central arrangements. 
Fix $H_1\in\A_1$ and $H_2\in\A_2$. Assume that 
\begin{itemize}
\item 
$\A_1$ is free, 
\item 
$(\A_1^{H_1}, m^{H_1})\simeq (\A_2^{H_2}, m^{H_2})$, and 
\item 
$b_2(\A_1)=b_2(\A_2)$. 
\end{itemize}
Then $\A_2$ is also free. 
\end{cor}

In Theorem \ref{thm:main} (3) and Corollary 
\ref{cor:main}, we can not drop the assumption 
that the multirestriction $(\A^{H_0}, m^{H_0})$ 
is free. Indeed, there exists non-free arrangement 
$\A$ such that $\chi(\dhA, t)=\chi((\A^{H_0}, m^{H_0}), t)$. 

\begin{example}
\label{ex:1}
Let $\A_1$ be a central arrangement in $\C^4$ 
defined by 
$$
x(x-w)y(y-w)(x+y+z)(x-y+z)zw. 
$$
Let $H_0=\{w=0\}$. Then $\dhA_1$ is an affine arrangement 
in $\C^3$ defined by 
$$
x(x-1)y(y-1)(x+y+z)(x-y+z)z
$$
whose characteristic polynomial is 
$\chi(\dhA_1, t)=t^3-7t^2+18t-17$. On the 
other hand, the multirestriction 
$(\A_1^{H_0}, m^{H_0})$ is defined by 
$x^2y^2(x+y+z)(x-y+z)z$. The characteristic 
polynomial is 
$\chi((\A_1^{H_0}, m^{H_0}),t)=t^3-7t^2+18t-17$ 
(\cite{ATW}) and 
we have $\sigma(\A_1^{H_0}, m^{H_0})=b_2(\dhA_1)=18$. 
However, since the characteristic polynomial does not 
factor, both $\A_1$ and $(\A_1^{H_0}, m^{H_0})$ are non-free. 
\end{example}

\begin{rem}
\label{rem:formal}
Recall that the arrangement $\A$ is called formal if every 
linear dependence $t_1\alpha_1+\dots+t_n\alpha_n=0$ of defining 
equations is a linear combination of 
three terms dependences 
$t_i\alpha_i+t_j\alpha_j+t_k\alpha_k=0$, in other 
words, linear dependences are generated by codimension 
two flats $L_2(\A)$. In \cite{yuz}, Yuzvinsky proved that 
free arrangements are formal. 
Theorem \ref{thm:main} shows that the freeness of $\A$ is 
characterized by combinatorial structures in codimension two 
and a multirestriction. Our results seem to have some 
relations with formality. However it is not clear yet. 
\end{rem}

\section{Related results}
\label{sec:related}

In general, $\chi(\dhA, t)$ and $\chi((\A^{H_0}, m^{H_0}), t)$ 
are not equal. However, under some assumptions on locally freeness, 
they are almost equal (they are equal except for 
the constant terms). 
We will give two different proofs for the following result. 

\begin{theorem}
If $\A$ is locally free along $H_0$, then 
$\chi(\dhA, t)-\chi((\A^{H_0}, m^{H_0}), t)\in\Z$. 
\end{theorem}

\subsection{First proof}

Let $1\leq k\leq \ell-2$. We shall prove 
$b_k(\dhA)=\sigma_k(\A^{H_0}, m^{H_0})$. 
From (\ref{eq:betti}), we have 
\begin{equation}
b_k(\dhA)=
\sum_{X\in L_k(\A^{H_0})}
\left(
\sum_{Y\in\rho^{-1}(X)}
b_k((\dhA)_Y)
\right). 
\end{equation}
Since $\A$ is locally free along $H_0$, 
$\A_X$ is free. Hence 
\begin{equation}
\sigma_k(\A_X^{H_0}, m_X^{H_0})=
b_k((\dhA_X))
=
\sum_{Y\in\rho^{-1}(X)}b_k((\dhA)_Y). 
\end{equation}
From local-global formula, we have 
$b_k(\dhA)=\sigma_k(\A^{H_0}, m^{H_0})$. \owari

\subsection{Second proof}

We first recall restriction maps for logarithmic forms 
following \cite{Sc, Y2}. 
Let us fix coordinates 
$z_1, \dots, z_\ell$ of $V$ so that 
$H_0=\{z_\ell=0\}$ (as in \S\ref{sec:simple}). 

A logarithmic differential form $\omega\in\Omega^p(\A)$ 
can be expressed as 
$$
\omega=
\omega_1+\frac{dz_\ell}{z_\ell}\wedge\omega_2, 
$$
where $\omega_1, \omega_2$ are rational 
differential forms generated by $dz_1, \dots, dz_{\ell-1}$. 
We can define the restriction map 
$\res_{H_0}^p:\Omega^p(\A)
\longrightarrow
\Omega^p(\A^{H_0}, m^{H_0})$
by 
$$
\omega_1+\frac{dz_\ell}{z_\ell}\wedge\omega_2
\longmapsto
\omega_1|_{H_0}. 
$$
The image of the map $\res_{H_0}^p$ is denoted by 
$
\res_{H_0}^p(\Omega^p(\A))=
M^p\subset
\Omega^p(\A^{H_0}, m^{H_0}), 
$
and its cokernel by $C^p$. We have the 
exact sequence 
$$
0
\longrightarrow 
M^p
\longrightarrow 
\Omega^p(\A^{H_0}, m^{H_0})
\longrightarrow 
C^p
\longrightarrow 
0. 
$$
\begin{prop}[\cite{Y2}]
If $\A$ is free, then $\res_\Hzero^p$ is surjective. 
\end{prop}
Define 
$\Phi(C^\bullet; x, y)$ to be 
$$
\Phi(C^\bullet; x, y)
=
\sum_{p=0}^{\ell-1}
P(C^p, x)y^p. 
$$
Then we have (see \cite{Sc})
\begin{equation}
\begin{split}
\chi((\A^{H_0}, m^{H_0}); t)-
\chi_0(\A, t)&=
\lim_{x\rightarrow 1}
\Phi(C^\bullet; x, t(1-x)-1)\\
&=
\lim_{x\rightarrow 1}
\sum_{p=0}^{\ell-1}
P(C^p, x)(t(1-x)-1)^p\\
&=
\lim_{x\rightarrow 1}
\sum_{p=0}^{\ell-1}
P(C^p, x)
\sum_{k=0}^p
(-1)^{p-k}
\begin{pmatrix}
p\\
k
\end{pmatrix}
t^k(1-x)^k. 
\\
&=
\lim_{x\rightarrow 1}
\sum_{k=0}^{\ell-1}
\begin{pmatrix}
p\\
k
\end{pmatrix}
t^k(1-x)^k
\left(
\sum_{p\geq k}
(-1)^{p-k}
P(C^p, x)
\right).
\end{split}
\end{equation}

Now we assume that $\A$ is locally free along $H_0$. 
Then the cokernel of the restriction map $C^p$ is 
supported on the origin $0\in H_0$. Therefore 
$\dim_\K C^p$ is finite dimensional, and the 
Hilbert series $P(C^p, x)$ is a (Laurent) polynomial. 
Hence $\lim_{x\rightarrow 1}(1-x)^kP(C^p, x)=0$ if 
$k\geq 1$. Thus we have 
$$
\chi((\A^{H_0}, m^{H_0}); t)-
\chi_0(\A, t)
=
\sum_{p=0}^{\ell-1}
(-1)^{p}\dim_\K C^p. 
$$
\owari

\begin{cor}
Let $\ell\geq 4$. Assume that 
the multirestriction $(\A^{H_0}, m^{H_0})$ is locally free, 
i.e., for any $0\neq X\subset H_0$, $(\A_X^{H_0}, m_X^{H_0})$ is 
free. Then the following are equivalent. 
\begin{itemize}
\item[(i)] $\A$ is locally free along $H_0$. 
\item[(ii)] $\chi(\dhA, t)-\chi((\A^{H_0}, m^{H_0}), t)\in\Z$. 
\item[(iii)] $b_2(\dhA)=\sigma_2(\A^{H_0}, m^{H_0})$. 
\end{itemize}
\end{cor}
The proof is similar to that of Theorem \ref{thm:main}. 

\begin{example}
\label{ex:2}
Let $\A_2$ be a central arrangement in $\C^4$ 
defined by 
$$
x(x-w)y(y-w)(x+y+z)(x-y+z)(z-w)w. 
$$
(Note that the 7-th hyperplane is different from 
$\A_1$ in Example \ref{ex:1}.) 
Let $H_0=\{w=0\}$. Then $\dhA_2$ is an affine arrangement 
in $\C^3$ defined by 
$$
x(x-1)y(y-1)(x+y+z)(x-y+z)(z-1)
$$
whose characteristic polynomial is 
$\chi(\dhA_2, t)=t^3-7t^2+18t-19$. On the 
other hand, the multirestriction 
$(\A_1^{H_0}, m^{H_0})$ is defined by 
$x^2y^2(x+y+z)(x-y+z)z$, which is the same 
arrangement with Example \ref{ex:1}. 
The characteristic polynomial is 
$\chi((\A_2^{H_0}, m^{H_0}),t)=t^3-7t^2+18t-17$. 
Since $(\A_2^{H_0}, m^{H_0})$ is rank three, hence 
locally free, and 
$\sigma(\A_2^{H_0}, m^{H_0})=b_2(\dhA_2)=18$, 
$\A_2$ is locally free along $H_0$. 
\end{example}

\vspace{5mm}

\noindent
Takuro Abe\\
Department of Mechanical Engineering and Science, \\
Kyoto University, \\
Yoshida Honmachi, Sakyo-ku, Kyoto 6068501, Japan. \\
abe.takuro.4c@kyoto-u.ac.jp



\bigskip

\noindent
Masahiko Yoshinaga\\
Department of Mathematics, \\
Kyoto University, \\
Kitashirakawa Oiwakecho, Sakyo-ku, Kyoto, 606-8502, Japan\\
mhyo@math.kyoto-u.ac.jp

\end{document}